# PREFERENCES RELATIONS AND CONSUMER THEORY

**FERREIRA, Manuel Alberto M.** (manuel.ferreira@iscte.pt)

**FILIPE, José António** (jose.filipe@iscte.pt)

INSTITUTO UNIVERSITÁRIO DE LISBOA (ISCTE-IUL),
BUSINESS RESEARCH UNIT (BRU-IUL) AND INFORMATION SCIENCES,
TECHNOLOGIES AND ARCHITECTURE RESEARCH CENTER (ISTAR-IUL),
AV. DAS FORÇAS ARMADAS
1649-026 LISBOA, PORTUGAL

**Abstract.** In this paper an interesting application of mathematics in economics is presented: the formulation of the theory of consumer basic problem, grounded on the concept of preferences relation and operationalized with optimization tools.

**Keywords.** Preferences relation, optimization, consumer theory, demand.

*MSC 2010:* 49K20.

## 1    Introduction

This work target is to present a mathematical formulation of the consumer theory as rigorous as possible. This is done using mainly the concept of preferences relation and then optimization mathematical tools.
 So, in the next section, the preferences relation concept and the respective axioms are presented and the model to be used is outlined.
In Section 3 is established the consumer theory and some important consequent results are shown.
Through Section 4 it is made a function $u(\cdot)$, which plays a key role in this work, short study.
In the end there is a short section of conclusions and bibliography on these subjects.

## 2    Preferences Relations

To use mathematical tools to deal with the concept of consumer basket, consider a non-empty convex closed subset $X$ in $\mathbb{R}_+^n$ where is defined a preferences relation, symbolized $\succcurlyeq$, as follows:
    **Definition 2.1**
For any $x$, $y$ and $z$ belonging to $X$,

*i)* $x \succcurlyeq y \ \vee \ y \succcurlyeq x$ (Completeness),

*ii)* $x \succcurlyeq y \ \wedge \ y \succcurlyeq z \implies x \succcurlyeq z$ (Transitivity)

*iii)* For any $y \in X$, $\{x \in X : x \succcurlyeq y\}$ and $\{x \in X : y \succcurlyeq x\}$ are closed subsets (Continuity),

*iv)* $x \geq y$ (that is: $x_i \geq y_i, i = 1, 2, \ldots n$) and $x \neq y \implies x \succ y$ (Strong Monotonicity).■

**Note**:

- $x \succ y$ if $x \succcurlyeq y$ and $y \succcurlyeq x$ is false,

- In economics language $X$ is the set of a consumer possible baskets and $\succcurlyeq$ represents that consumer preferences lengthwise this set of baskets.

**Definition 2.2**

$x, y \in X$ are indifferent, and it is denoted $x \sim y$, when $x \succcurlyeq y$ and $y \succcurlyeq x$.■

**Definition 2.3**[1]

A basket $x \in X$ is said to be redundant when there is $y \in X$ such that $x \geq y, y \neq x$, and $x \sim y$. ■

To present some results, to the Definition 1.1 must be added the following two axioms:

*v)* $x, y \in X, x \neq y$ and $x, y$ not redundant such that $x \sim y \implies (1 - t)x + ty \succ x, \forall t \in [0,1]$.■

**Note:**

-That is: the convex linear combination of two different, not redundant, and indifferent baskets is preferred to one of them isolated.

*vi)* $\forall x \in X, \forall \varepsilon > 0, \exists y \in X$ such that $\|x - y\| < \varepsilon$ and $y \succ x$.■

**Note**:

-That is: in any neighborhood of a basket there is at least one other basket preferable to it.

**Theorem 2.1**

Supposing that $\succcurlyeq$ fulfills the axioms *i), ii), iii)* and *iv)* in Definition 2.1, if $X = \mathbb{R}_+{}^n$ there is a continuous function $u : X \to \mathbb{R}$ such that $u(x) \geq u(y)$ if and only if $x \succcurlyeq y, \forall x, y \in X$.

**Dem**:

Call **1** the vector with the whole $n$ coordinates equal to 1. For $x \in X$ define $A = \{t \in \mathbb{R}_+ : t\mathbf{1} \succcurlyeq x\}$ and $B = \{t \in \mathbb{R}_+ : x \succcurlyeq t\mathbf{1}\}$. According to *iv)* $A$ and $B$ are nonempty and by *iii)* are both closed since the function $\varphi : \begin{array}{c} \mathbb{R}_+ \to \mathbb{R}_+{}^n \\ t \to t\mathbf{1} \end{array}$ is continuous and, here, $A = \varphi^{-1}(\{y \in X : y \succcurlyeq x\})$ and $B = \varphi^{-1}(\{y \in X : x \succcurlyeq y\})$. After *i)* it results $A \cup B = \mathbb{R}_+$ and by the $\mathbb{R}_+$ connexity [2] there is $t \in \mathbb{R}$ fulfilling $t\mathbf{1} \sim x$. Again, by *iv)* this $t$ is unique. So define $u : X \to \mathbb{R}$ as $u(x) = t$. Besides $u^{-1}([t_0, \infty]) = \{x \in X : u(x) \geq t_0\} = \{x \in X : x \succcurlyeq u^{-1}(t_0)\}$ and $u^{-1}([0, t_0]) = \{x \in X : u^{-1}(t_0) \succcurlyeq x\}$ are closed and then it results that $u$ is continuous[3].■

**Note**:

-More generally it is also true that "If $X \subset \mathbb{R}_+$ is a connex set and $\succcurlyeq$ is a preferences relation in $X$ fulfilling *i), ii)* and *iii)*, there is a continuous function $u : X \to \mathbb{R}$ such that $u(x) \geq u(y) \iff x \succcurlyeq y, \forall x, y \in X$", see Debreu (1959),

---

[1] See Simonsen (1989).

[2] A set is connex if it cannot be represented as the union of two separable sets both non-empty.

[3] Because in this case the whole closed sets in $[0, \infty)$ will have inverse image closed. For details see (Cysne and Moreira, 1997), page 95.

-This theorem shows that under axioms *i), ii), iii)* and *iv)* it is possible to determine a numerical scale for the consumer preferences what is extremely important in Decision Theory.

## 3    Consumer Theory

Begin by a formulation, among many equivalents, of the theory of consumer basic problem:

-Call $r$ the consumer income and $p = (p_1, p_2, \ldots, p_n)$ the $n$ goods vector prices. If $u: X \to \mathbb{R}$ represents the consumer preferences, the preferences maximization problem can be mathematically formulated as

$$\max \; u(x)$$
$$s.to \; \langle p, x \rangle \leq r \quad (3.1).$$

It is supposed that

$$A = \{x \in X : \langle p, x \rangle \leq r\},$$

the opportunities set, is nonempty. ■

**Note**:

-Being $u$ continuous and $p \in \mathbb{R}_+^{*\,n}$, so being A compact, the Weierstrass theorem guaranties that (3.1) is possible.

**Definition 3.1**

The indirect utility function $v(p, r)$ is defined through the problem:

$$v(p, r) = \max \; u(x)$$
$$s.to \; \langle p, x \rangle \leq r, x \in X \quad (3.2). ■$$

**Note**:

-It is easy to check that if $x$ is a solution for (3.2) $\langle p, x \rangle = r$ since either *vi)* or *iv)* are satisfied for $X = \mathbb{R}_+^{\,n}$,

- $x$ is not redundant,

- If, in addition, *v)* is satisfied $x$ is unique.

**Definition 3.2**

The vector only solution $x$ for the problem (3.2), given $p$ and r: $x(p, r)$ is the Marshallian demand. ■

**Note**:

-$x_i(p, r)$ is the Marshallian demand for the good $i, i = 1, 2, \ldots, n$.

**Theorem 3.1**

If consumer preferences satisfy *i), ii), iii)* and *vi)*, the Marshallian demand $x_i : \mathbb{R}_+^{*\,n} \times \mathbb{R}_+ \to \mathbb{R}$ is continuous $\forall i = 1, 2, \ldots, n$.

**Dem**: See Simonsen (1989). ■

## 4    The Function **u**$(\cdot)$

Consider the problem (3.1) following version:

$$\max \; u(x)$$
$$s.to \; \langle p, x \rangle = r \quad (4.1).$$

It is a constrained optimization problem. The Lagrange's multipliers method may be applied in its solution, see for instance (Ferreira and Amaral, 2002). The Lagrangean function is

$$\mathcal{L}(x_1, x_2, \ldots, x_n, \lambda) = u(x_1, x_2, \ldots, x_n) + \lambda(r - x_1 p_1 - x_2 p_2 - \cdots - x_n p_n) \qquad (4.2).$$

The first order conditions are

$$\frac{\partial u}{\partial x_1} - \lambda p_1 = 0$$

$$\frac{\partial u}{\partial x_2} - \lambda p_2 = 0$$

$$\vdots$$

$$\frac{\partial u}{\partial x_n} - \lambda p_n = 0$$

$$x_1 p_1 + x_2 p_2 + \cdots + x_n p_n = r.$$

Using the first and the last conditions it is obtained, $x_1 \neq 0$,

$$\frac{\partial u}{\partial x_1} - \lambda \left( \frac{Y}{x_1} - \frac{x_2}{x_1} p_2 - \cdots - \frac{x_n}{x_1} p_n \right) = 0$$

or

$$x_1 \frac{\partial u}{\partial x_1} + \lambda x_2 p_2 + \cdots + \lambda x_n p_n = \lambda r.$$

As $\lambda p_i = \frac{\partial u}{\partial x_i}, i = 1, 2, \ldots, n$, the former equation becomes

$$x_1 \frac{\partial u}{\partial x_1} + x_2 \frac{\partial u}{\partial x_2} + \cdots + x_n \frac{\partial u}{\partial x_n} = \lambda r$$

and noting that

$$\lambda = \frac{1}{p_1} \frac{\partial u}{\partial x_1}$$

it assumes the form

$$\left( x_1 - \frac{Y}{p_1} \right) \frac{\partial u}{\partial x_1} + x_2 \frac{\partial u}{\partial x_2} + \cdots + x_n \frac{\partial u}{\partial x_n} = 0 \qquad (4.3).$$

It is a first order homogeneous partial derivatives equation, (Ferreira and Amaral, 2005). Solving it:

$$u(x_1, x_2, \ldots, x_n) = F\left(\frac{r}{p_1 + p_2 + \cdots + p_n}\right) \qquad (4.4).$$

**Note:**

-If $u_1$ and $u_2$ solve (4.3), $\alpha_1 u_1 + \alpha_2 u_2, \alpha_1, \alpha_2 \in \mathbb{R}$ also solve it:

$\left(x_1 - \frac{r}{p_1}\right)\left(\alpha_1 \frac{\partial u_1}{\partial x_1} + \alpha_2 \frac{\partial u_2}{\partial x_1}\right) + x_2\left(\alpha_1 \frac{\partial u_1}{\partial x_2} + \alpha_2 \frac{\partial u_2}{\partial x_2}\right) + \cdots + x_n\left(\alpha_1 \frac{\partial u_1}{\partial x_n} + \alpha_2 \frac{\partial u_2}{\partial x_n}\right) = \alpha_1\left(\left(x_1 - \frac{r}{P_1}\right)\frac{\partial u_1}{\partial x_1} + x_2 \frac{\partial u_1}{\partial x_2} + \cdots + x_n \frac{\partial u_1}{\partial x_n}\right) + \alpha_2\left(\left(x_1 - \frac{r}{P_{x_1}}\right)\frac{\partial u_2}{\partial x_1} + x_2 \frac{\partial u_2}{\partial x_2} + \cdots + x_n \frac{\partial u_2}{\partial x_n}\right) = \alpha_1 0 + \alpha_2 0 = 0$. So, the set of functions is linear and, of course, a linear combination of solutions is also a solution,

-$F(.)$ is any differentiable function,

-$r = x_1 P_{x_1} + x_2 P_{x_2} + \cdots + x_n P_{x_n}$,

-It is easy to check that (4.4) is a solution for (4.3) substituting directly,

-The expression (4.4) evidences the utility functional dependency from the whole goods and the income.∎

**Some Examples:**

-One concretization of (4.4) may be

$$u(x_1, x_2, \ldots, x_n) = \alpha^{\frac{r}{p_1 + p_2 + \cdots + p_n}} \quad (4.5)$$

for which $U(1, 1, \ldots, 1) = \alpha$. That is $\alpha$ is a standard utility: the value of the utility when unitary quantities of every good are used.∎

-Defining $z_i = \alpha^{x_i}, i = 1, 2, \ldots, n$, (4.5) becomes

$$u(x_1, x_2, \ldots, x_n) = z_1{}^{\alpha_1} z_2{}^{\alpha_2} \ldots z_n{}^{\alpha_n} \quad (4.6)$$

where $\alpha_i = \frac{p_i}{p}, i = 1, 2, \ldots, n$ and $p = \sum_{i=1}^{n} p_i$ are the standardized prices of each good. Formally, (4.6) is a Cobb-Douglas function. That is: in terms of the standard utility base exponential of each quantity, (4.5) assumes the Cobb-Douglas utility form.∎

-Another example is

$$u(x_1, x_2, \ldots, x_n) = \beta \frac{r}{p_1 + p_2 + \cdots + p_n} \qquad (4.7)$$

being now $\beta$ the standard utility. ∎

-The expression (4.7) may be written as

$$u(x_1, x_2, \ldots, x_n) = \beta \sum_{i=1}^{n} \alpha_i x_i \quad (4.8)$$

with $\alpha_i$ defined above. So, in this case, the utility is given by a linear combination of the quantities, which weights are the standardized prices, multiplied by the standard utility. ∎

-And, finally, as a last example

$$u(x_1, x_2, \ldots, x_n) = \gamma \, ln \frac{r}{p_1 + p_2 + \cdots + p_n} \quad (4.9)$$

for which $U(1, 1, \ldots, 1) = 0$ and $U(e, e, \ldots, e) = \gamma$ .

Evidently

$$u(x_1, x_2, \ldots, x_n) = \, ln(\sum_{i=1}^{n} \alpha_i x_i)^\gamma \quad (4.10). ∎$$

## 5    Concluding Remarks

In this paper, economic concepts were presented through a rigorous formulation, allowing that some economic results were demonstrated in the mathematical sense. Plays here a fundamental role the concept of preferences relation that is axiomatically formulated.

## Bibliography

[1] BORCH, K. (1974), *The Mathematical Theory of Insurance*. Lexington, Mass.: Lexington Books, D.C. Heath and Co..

[2] BORCH; K. (1990), *Economics of Insurance*. Edited by Aase K. K. and Sandmo A., Amsterdam, North-Holland.

[3] CHAVAGLIA, J., FILIPE, J. A., FERREIRA, M. A. M. (2016), "Neuroeconomics and reinformance learning an exploratory analysis", *APLIMAT 2016-15th Conference on Applied Mathematics 2016, Proceedings 2016*, Bratislava, 527-534.

[4] CYSNE, R. P., MOREIRA, H. A. (1997), *Curso de Matemática para Economistas*. São Paulo, Editora Atlas, S. A..

[5] DEBREU,G. (1959), *Theory of value and axiomatic analysis of economic equilibrium*. Monograph, 17, Cowles Foundation.

[6] FERREIRA, M. A. M. (1991), "Equações com derivadas parciais e funções de utilidade", *Revista de Gestão,* XI-XII, 63.

[7] FERREIRA, M. A. M. (2015), "The Hahn-Banach theorem in convex programming", *APLIMAT 2015-14th Conference on Applied Mathematics 2015, Proceedings 2015,* Bratislava, 283-291.

[8] FERREIRA, M. A. M. (2015a), "The minimax theorem as Hahn-Banach theorem consequence", *Acta Scientiae et Intellectus*, 1 (1), 58-65.

[9] FERREIRA, M. A. M., AMARAL, I. (2002), *Cálculo Diferencial em $\mathbb{R}^n$* . Edições Sílabo-5ª Edição, Lisboa.

[10]     FERREIRA, M. A. M., AMARAL, I. (2005), *Integrais Múltiplos e Equações Diferenciais*. Edições Sílabo-5ª Edição, Lisboa.

[11]     FERREIRA, M. A. M., ANDRADE, M. (2011), "A note on partial derivatives equations and utility functions", *Journal of Economics and Engineering*, 2(1), 23-24.

[12]     FERREIRA, M. A. M., FILIPE, J. A. (2014), "Convex sets strict separation in Hilbert spaces", *Applied Mathematical Sciences*, (61-64), 3155-3160. **DOI**: 10.12988/ams.2014.44257

[13]     FERREIRA, M. A. M., FILIPE, J. A. (2015), "A ratio to evaluate harvest procedures management in an economic system where resources dynamics is ruled by an ornstein-uhlenbeck process", *Applied Mathematical Sciences*, 9(25-28), 1213-1219. **DOI**: 10.12988/ams.2015. 410871

[14]     FERREIRA, M. A. M., MATOS, M. C. (2014), "Convex sets strict separation in the minimax theorem", *Applied Mathematical Sciences*, 8(33-36), 1781-1787. **DOI**: 10.12988/ams.2014.4271

[15]     FERREIRA, M. A. M., FILIPE, J. A., ANDRADE, M. (2013), "A note on partial derivatives equations and utility functions (Revisited)", *Journal of Economics and Engineering*, 2(1), 11-13.

[16]     FERREIRA, M. A. M., ANDRADE, M., MATOS, M. C., FILIPE, J. A., COELHO, M. (2012), "Kuhn-Tucker's Theorem-the Fundamental Result in Convex Programming Applied to Finance and Economic Sciences", *International Journal of Latest Trends in Finance & Economic Sciences*, 2 (2), 111-116.

[17]     FILIPE, J. A., FERREIRA, M. A. M. (2013), "Chaos in humanities and social sciences: An approach", *APLIMAT 2013-12th Conference on Applied Mathematics 2013, Proceedings 2013,* Bratislava.

[18]     NEUMÄRKER, B. (2007), "Neuroeconomics and the Economic Logic of Behavior", *Analyse & Kritik*, 29, 60-85.

[19]     ROBERTS, B., SCHULZE, D. L. (1973), *Modern Mathematics and Economic Analysis*. New York: W. W. Norton & Company, Inc..

**Current address**

**Manuel Alberto M. Ferreira, Professor Catedrático**
INSTITUTO UNIVERSITÁRIO DE LISBOA (ISCTE-IUL),
BUSINESS RESEARCH UNIT (BRU-IUL) AND INFORMATION SCIENCES, TECHNOLOGIES AND ARCHITECTURE RESEARCH CENTER (ISTAR-IUL), AV. DAS FORÇAS ARMADAS
1649-026 LISBOA, PORTUGAL
TELEFONE: + 351 21 790 37 03
FAX: + 351 21 790 39 41
E-MAIL: manuel.ferreira@iscte.pt

**José António Filipe, Professor Auxiliar**
INSTITUTO UNIVERSITÁRIO DE LISBOA (ISCTE-IUL),
BUSINESS RESEARCH UNIT (BRU-IUL) AND INFORMATION SCIENCES, TECHNOLOGIES AND ARCHITECTURE RESEARCH CENTER (ISTAR-IUL), AV. DAS FORÇAS ARMADAS
1649-026 LISBOA, PORTUGAL
TELEFONE: + 351 21 790 34 08
FAX: + 351 21 790 39 41
E-MAIL: jose.filipe@iscte.pt